\documentclass[12pt]{article}
\usepackage[]{amsmath,amssymb}
\usepackage{amscd}
\usepackage{latexsym}
\usepackage{cite}

\newtheorem{definition}{Definition}[section]
\newtheorem{theorem}[definition]{Theorem}
\newtheorem{lemma}[definition]{Lemma}

\newtheorem{problem}[definition]{Problem}
\newtheorem{note}[definition]{Note}

\def\C{\mathbb C}

\def\F{\mathbb F}
\def\K{\mathbb F}

\begin{document}
\title{\bf The Rahman polynomials and\\ 
the Lie algebra $\mathfrak{sl}_3(\C)$
}
\author{
Plamen Iliev\footnote{Supported in part by NSF grant DMS-0901092.}
\, 
and
Paul Terwilliger}
%\footnote{This author gratefully acknowledges 
%support from the FY2007 JSPS Invitation Fellowship Program
%for Reseach in Japan (Long-Term), grant L-07512.}
%}
\date{}
%to get date printout, comment out above line

\maketitle
\begin{abstract}
We interpret the Rahman polynomials
in terms of the Lie
algebra $\mathfrak{sl}_3(\C)$.
Using the parameters of the polynomials 
we define two Cartan subalgebras for
 $\mathfrak{sl}_3(\C)$, denoted $H$ and $\tilde H$.
We display an antiautomorphism $\dagger$ of
$\mathfrak{sl}_3(\C)$ that fixes each
element of $H$ and each element of $\tilde H$.
We consider a certain finite-dimensional irreducible
$\mathfrak{sl}_3(\C)$-module $V$ consisting
of homogeneous polynomials in three variables.
We display a nondegenerate symmetric bilinear form
$\langle \,,\,\rangle$ on $V$
such that
$\langle \beta \xi,\zeta\rangle = \langle \xi,\beta^\dagger \zeta\rangle
$
for all $ \beta \in
\mathfrak{sl}_3(\C)$ and $\xi,\zeta \in V$.
We display two bases for $V$;
one diagonalizes $H$ and the other diagonalizes
$\tilde H$.
Both bases are orthogonal with respect to
$\langle \,,\,\rangle$.
We show that when $\langle\,,\,\rangle$
is applied to a vector in
each basis, the result is 
a trivial factor
times a Rahman polynomial evaluated at an appropriate argument.
Thus for both transition matrices between the bases
each entry is described by a 
Rahman polynomial.
From these results
we recover
the previously known orthogonality relation for 
the Rahman polynomials.
We also obtain two seven-term recurrence relations
satisfied by the Rahman polynomials, along with 
the corresponding relations satisfied by the dual polynomials.
These recurrence relations show that the Rahman polynomials
are bispectral.
In our theory the roles of
$H$ and $\tilde H$
are interchangable,
and for us 
this explains the duality and bispectrality of the Rahman polynomials.
We view the action of $H$ and $\tilde H$ on $V$ as
a rank 2 generalization of a Leonard pair.

 \end{abstract}

\section{The Rahman polynomials}\label{se1}
\noindent We begin by recalling the Rahman
polynomials
\cite{grun},
\cite{rahman}.
In what follows $\lbrace p_i\rbrace_{i=1}^4$
denote complex numbers. They are essentially
arbitrary, although certain combinations are forbidden
in order to avoid dividing by zero.
Define
\begin{eqnarray*}
&&t =\frac{(p_1+p_2)(p_1 +p_3)}{p_1(p_1+p_2+p_3+p_4)},
\qquad
u =\frac{(p_1+p_3)(p_3 +p_4)}{p_3(p_1+p_2+p_3+p_4)},
\\
&&v =\frac{(p_1+p_2)(p_2 +p_4)}{p_2(p_1+p_2+p_3+p_4)},
\qquad
w =\frac{(p_2+p_4)(p_3 +p_4)}{p_4(p_1+p_2+p_3+p_4)}.
\end{eqnarray*}
%It turns out that
%$(1-t^{-1})(1-w^{-1})$ and $(1-u^{-1})(1-v^{-1})$
%are equal; denote the common value by $\nu^{-1}$.
Fix an integer $N\geq 0$ and
let $a,b,c,d$ denote
mutually commuting
indeterminates.
Define 
 \begin{eqnarray*}
P (a,b,c,d)
=
\sum_{\genfrac{}{}{0pt}{} {0\leq i,j,k,\ell}{i+j+k+\ell\leq N}}
\frac{
(-a)_{i+j}
(-b)_{k+\ell}
(-c)_{i+k}
(-d)_{j+\ell}
}
{i!j!k!\ell ! (-N)_{i+j+k+\ell}}
t^iu^j v^k w^\ell.
\end{eqnarray*}
We are using the shifted factorial notation
\begin{eqnarray*}
    (\alpha)_n = \alpha (\alpha +1)\cdots (\alpha+n-1)
    \qquad \qquad n=0,1,2,\ldots
\end{eqnarray*}
We will use the fact that for nonnegative integers 
$m,n$
\begin{eqnarray}
(-m)_n = 0 \quad {\rm if} \quad n > m.
\label{eq:principle-intro}
\end{eqnarray}
For nonnegative integers $m,n$ whose sum is at most $N$
the corresponding Rahman polynomial is
$P (m,n,c,d)$ in the variables $c,d$,
and the corresponding dual Rahman polynomial 
is
$ P (a,b,m,n)$ in the variables $a,b$
\cite[Section 2]{rahman}.
The Rahman polynomials and their duals satisfy
an orthogonality relation which we now describe.
Define
\begin{eqnarray*}
\nu 
&=& \frac
{(p_1+p_2)(p_1+p_3)(p_2+p_4)(p_3+p_4)}
{(p_1 p_4-p_2 p_3)^2}.
\end{eqnarray*}
Define  $\eta_0 = \nu^{-1}$ and
\begin{eqnarray*}
\eta_1 &=& \frac{p_1 p_2 (p_1+p_2+p_3+p_4)}{(p_1+p_2)(p_1+p_3)(p_2+p_4)},
\\
\eta_2 &=& \frac{p_3 p_4 (p_1+p_2+p_3+p_4)}{(p_1+p_3)(p_2+p_4)(p_3+p_4)}.
\end{eqnarray*}
Define
$\tilde \eta_0 = \nu^{-1}$ and
\begin{eqnarray*}
\tilde \eta_1 &=& \frac{p_1 p_3 (p_1+p_2+p_3+p_4)}{(p_1+p_2)(p_1+p_3)(p_3+p_4)},
\\
\tilde \eta_2 &=& \frac{p_2 p_4 (p_1+p_2+p_3+p_4)}{(p_1+p_2)(p_2+p_4)(p_3+p_4)}.
\end{eqnarray*}
A short computation shows
\begin{eqnarray*}
 \eta_0+ \eta_1+ \eta_2=1,
 \qquad \qquad 
\tilde \eta_0+ \tilde \eta_1+ \tilde \eta_2=1.
\end{eqnarray*}
For notational convenience 
define $k_i=\nu \tilde \eta_i$
and $\tilde k_i=\nu \eta_i$ for
$0 \leq i \leq 2$, so that $k_0=1$ and
 $\tilde k_0=1$.

\medskip
\noindent We credit the following result to Rahman and Hoare \cite{rahman};
however a similar and more general theorem
was given earlier in
\cite[Theorem~1.1]{tanaka}.

\begin{theorem}
\label{thm:rahmanorthog}
{\rm (Rahman and Hoare\cite{rahman})}
Fix nonnegative integers $s$, $t$ whose sum is at most $N$,
and nonnegative integers $\sigma$, $ \tau$ whose sum is at most $N$.
Then both 
\begin{eqnarray}
\label{eq:orthog1-intro}
\sum_{\genfrac{}{}{0pt}{} {0\leq i,j,k}{i+j+k=N}}
 P (j,k,s,t)
 P (j,k,\sigma,\tau)
\tilde\eta^i_0
\tilde \eta^j_1
\tilde \eta^k_2
\binom{N}{i \; j\; k} 
&=&
\frac{\delta_{s \sigma} \delta_{t \tau}
}
{
\tilde k^{s}_1 \tilde k^{t}_2
}
\binom{N}{r \; s\; t}^{-1},
\\
\label{eq:orthog2-intro}
\sum_{\genfrac{}{}{0pt}{} {0\leq i,j,k}{i+j+k=N}}
 P (s,t,j,k)
 P (\sigma,\tau,j,k)
\eta^i_0
 \eta^j_1
 \eta^k_2
\binom{N}{i \; j\; k} 
&=& \frac{\delta_{s \sigma} \delta_{t \tau}}
{
 k^{s}_1  k^{t}_2} 
\binom{N}{r \; s\; t}^{-1},
\end{eqnarray}
 where $r=N-s-t$.
\end{theorem}

\noindent
In this paper 
we interpret the Rahman polynomials 
in terms of the Lie
algebra $\mathfrak{sl}_3(\C)$.
Our results are summarized as follows.
Using 
the parameters $\lbrace p_i\rbrace_{i=1}^4$
we define two Cartan subalgebras for
$\mathfrak{sl}_3(\C)$, denoted $H$ and $\tilde H$.
We display an antiautomorphism $\dagger$ of
$\mathfrak{sl}_3(\C)$ that fixes each
element of $H$ and each element of $\tilde H$.
We consider an irreducible
$\mathfrak{sl}_3(\C)$-module $V$ consisting
of the homogeneous polynomials in three variables that
have total degree $N$.
We display a nondegenerate symmetric bilinear form
$\langle \,,\,\rangle$ on $V$
such that
$\langle \beta \xi,\zeta\rangle = \langle \xi,\beta^\dagger \zeta\rangle
$
for all $ \beta \in
\mathfrak{sl}_3(\C)$ and $\xi,\zeta \in V$.
We display two bases for $V$;
one diagonalizes $H$ and the other diagonalizes
$\tilde H$.
Both bases are orthogonal with respect to
$\langle \,,\,\rangle$.
We show that when $\langle\,,\,\rangle$
is applied to a vector in
each basis, the result is 
a trivial factor
times a Rahman polynomial evaluated at an appropriate argument.
Thus for both transition matrices between the bases
each entry is described by a 
Rahman polynomial.
%up to a trivial factor.
From these results we
obtain
an elementary proof of
Theorem 
\ref{thm:rahmanorthog}.
We also obtain
two seven-term recurrence relations
satisfied by the Rahman polynomials, along with the corresponding
relations satisfied by their duals.
These recurrence relations show that the Rahman polynomials
are bispectral, a feature hinted at by
Gr\"unbaum 
\cite{grun}.
We view the actions of $H$ and $\tilde H$ on
$V$ as a rank 2 generalization of 
a Leonard pair
\cite{LS99}; this is discussed along
with some open problems at the
end of the paper.

\section{The Lie algebra $\mathfrak{sl}_3(\C)$}\label{se2}
We will be discussing the Lie algebra
$\mathfrak{sl}_3(\C)$; for background information
see
\cite{carter}.
Let ${\rm {Mat}}_3(\C)$ denote the
$\C$-algebra consisting of the
$3 \times 3$ matrices that have all entries in 
$\C$. 
We index the rows and columns by $0,1,2$.
For $0 \leq i,j\leq 2$ let $e_{ij} $
denote the matrix in 
 ${\rm {Mat}}_3(\C)$
that has $(i,j)$-entry 1 and
all other entries 0.
The Lie algebra
$\mathfrak{sl}_3(\C)$ is the set of
matrices in
${\rm {Mat}}_3(\C)$ 
that have
trace 0, together with the Lie bracket
$\lbrack \beta,\gamma\rbrack=\beta \gamma -\gamma \beta$.
We will consider two Cartan subalgebras
of 
$\mathfrak{sl}_3(\C)$, denoted
$H$ and $\tilde H$.
The subalgebra $H$
consists of the diagonal matrices in
$\mathfrak{sl}_3(\C)$.
Define
\begin{eqnarray*}
\varphi= {\rm {diag}}(-1/3, 2/3, -1/3),
\qquad \qquad 
\phi  ={\rm {diag}}(-1/3, -1/3, 2/3).
\end{eqnarray*}
Then
$\varphi$, $\phi$ form a basis for $H$.
We will describe $\tilde H$ shortly.
Define
\begin{eqnarray*}
U=
\left(
\begin{array}{c c c }
1 & 
1 & 
1 
  \\
1 & 
1-t
&
 1-v 
  \\
1& 
1-u
&
1-w
\\
\end{array}
\right)
\end{eqnarray*}
and
\begin{eqnarray*}
W= {\rm {diag}}(\eta_0, \eta_1, \eta_2),
\qquad \qquad 
\tilde W  ={\rm {diag}}(\tilde \eta_0, \tilde \eta_1, \tilde\eta_2).
\end{eqnarray*}
A brief calculation shows
$\nu WU \tilde W U^t =  I$, 
where $I$ denotes the identity matrix and $U^t$ denotes the
transpose of $U$.
Define
\begin{eqnarray*}
\vartheta =
\frac
{(p_1+p_3)(p_2+p_4)}
{p_2 p_3-p_1 p_4},
\qquad \qquad 
\tilde \vartheta = \frac
{(p_1+p_2)(p_3+p_4)}
{p_2 p_3-p_1 p_4}
\end{eqnarray*}
and note that
$\vartheta \tilde \vartheta  = \nu $.
Define
$R = \tilde \vartheta \tilde W U^t$
and note that
$R^{-1} = \vartheta  W U$.
We have
\begin{eqnarray*}
R&=&
\left(
\begin{array}{c c c }
\frac{p_2p_3-p_1p_4}{(p_1+p_3)(p_2+p_4)} & 
\frac{p_2p_3-p_1p_4}{(p_1+p_3)(p_2+p_4)} & 
\frac{p_2p_3-p_1p_4}{(p_1+p_3)(p_2+p_4)} 
  \\
\frac{p_1 p_3 (p_1+p_2+p_3+p_4)}{(p_1+p_3)(p_2p_3-p_1p_4)}
& \frac{-p_3}{p_1+p_3}&
\frac{p_1}{p_1+p_3}  \\
\frac{p_2 p_4 (p_1+p_2+p_3+p_4)}{(p_2+p_4)(p_2p_3-p_1p_4)}
& \frac{p_4}{p_2+p_4}&
\frac{-p_2}{p_2+p_4}  \\
\end{array}
\right)
\end{eqnarray*}
and
\begin{eqnarray*}
R^{-1} &=&
\left(
\begin{array}{c c c }
\frac{p_2p_3-p_1p_4}{(p_1+p_2)(p_3+p_4)} & 
\frac{p_2p_3-p_1p_4}{(p_1+p_2)(p_3+p_4)} & 
\frac{p_2p_3-p_1p_4}{(p_1+p_2)(p_3+p_4)} 
  \\
\frac{p_1 p_2 (p_1+p_2+p_3+p_4)}{(p_1+p_2)(p_2p_3-p_1p_4)}
& \frac{-p_2}{p_1+p_2}&
\frac{p_1}{p_1+p_2}  \\
\frac{p_3 p_4 (p_1+p_2+p_3+p_4)}{(p_3+p_4)(p_2p_3-p_1p_4)}
& \frac{p_4}{p_3+p_4}&
\frac{-p_3}{p_3+p_4}  \\
\end{array}
\right).
\end{eqnarray*}
Note that 
$R W R^t = \tilde \vartheta \vartheta^{-1}\tilde W$.
For $0 \leq i,j\leq 2$ define
$\tilde e_{ij}= R e_{ij} R^{-1}$.
Define $\tilde H = R H R^{-1}$
and note that $\tilde H$ is a Cartan subalgebra of
$\mathfrak{sl}_3(\C)$.
Define
$\tilde \varphi= R \varphi R^{-1}$ and
$\tilde \phi= R \phi R^{-1}$.
Note that 
$\tilde \varphi$,
$\tilde \phi$ is a basis for
$\tilde H $.
We have
\begin{eqnarray*}
\tilde \varphi &=&
\frac{p_2(p_1p_4-p_2p_3)}{
(p_1+p_2)(p_1+p_3)(p_2+p_4)} e_{01}
\quad
+
\quad
\frac{p_1(p_2p_3-p_1p_4)}{
(p_1+p_2)(p_1+p_3)(p_2+p_4)} e_{02}
\\
&+&
\frac{p_1p_2p_3(p_1+p_2+p_3+p_4)}{(p_1+p_2)(p_1+p_3)(p_1p_4-p_2p_3)}
e_{10}
\quad
-
\quad
\frac{p_1 p_3}{
(p_1+p_2)(p_1+p_3)}
e_{12}
\\
&+&
\frac{p_1p_2p_4(p_1+p_2+p_3+p_4)}{(p_1+p_2)(p_2+p_4)(p_2p_3-p_1p_4)}
e_{20}
\quad
-
\quad
\frac{p_2 p_4}{
(p_1+p_2)(p_2+p_4)} 
e_{21}
\\
&+&
\biggl(\frac{p_2 p_3}{
(p_1+p_2)(p_1+p_3)}
-
\frac{p_1p_2(p_1+p_2+p_3+p_4)}{(p_1+p_2)(p_1+p_3)(p_2+p_4)}
\biggr) \varphi
\\
&+&
\biggl(\frac{p_1 p_4}{
(p_1+p_2)(p_2+p_4)}
-
\frac{p_1p_2(p_1+p_2+p_3+p_4)}{(p_1+p_2)(p_1+p_3)(p_2+p_4)}
\biggr) \phi,
\end{eqnarray*}
\begin{eqnarray*}
\tilde \phi &=&
\frac{p_4(p_2p_3-p_1p_4)}{
(p_1+p_3)(p_3+p_4)(p_2+p_4)}
e_{01}
\quad 
+ \quad
\frac{p_3(p_1p_4-p_2p_3)}{
(p_1+p_3)(p_3+p_4)(p_2+p_4)}
e_{02}
\\
&+&
\frac{p_1p_3p_4(p_1+p_2+p_3+p_4)}{(p_1+p_3)(p_3+p_4)(p_2p_3-p_1p_4)}
e_{10}
\quad
- 
\quad
\frac{p_1 p_3}{
(p_1+p_3)(p_3+p_4)}
e_{12}
\\
&+&
\frac{p_2p_3p_4(p_1+p_2+p_3+p_4)}{(p_2+p_4)(p_3+p_4)(p_1p_4-p_2p_3)}
e_{20}
\quad
-
\quad
\frac{p_2 p_4}{
(p_2+p_4)(p_3+p_4)}
e_{21}
\\
&+&
\biggl(
\frac{p_1 p_4}{
(p_1+p_3)(p_3+p_4)}
-
\frac{p_3p_4(p_1+p_2+p_3+p_4)}{(p_1+p_3)(p_3+p_4)(p_2+p_4)}
\biggr)
\varphi
\\
&+&
\biggl(
\frac{p_2 p_3}{
(p_2+p_4)(p_3+p_4)}
-
\frac{p_3p_4(p_1+p_2+p_3+p_4)}{(p_1+p_3)(p_3+p_4)(p_2+p_4)}
\biggr) \phi,
\end{eqnarray*}
\begin{eqnarray*}
 \varphi  &=&
\frac{p_3(p_1p_4-p_2p_3)}{
(p_1+p_2)(p_1+p_3)(p_3+p_4)}
\tilde e_{01}
\quad 
+ \quad
\frac{p_1(p_2p_3-p_1p_4)}{
(p_1+p_2)(p_1+p_3)(p_3+p_4)}
\tilde e_{02}
\\
&+&\frac{p_1p_2p_3(p_1+p_2+p_3+p_4)}{(p_1+p_2)(p_1+p_3)(p_1p_4-p_2p_3)}
\tilde e_{10}
\quad -
\quad
\frac{p_1 p_2}{
(p_1+p_2)(p_1+p_3)}
\tilde e_{12}
\\
&+& \frac{p_1p_3p_4(p_1+p_2+p_3+p_4)}{(p_1+p_3)(p_3+p_4)(p_2p_3-p_1p_4)}
\tilde e_{20}
\quad - \quad
\frac{p_3 p_4}{
(p_1+p_3)(p_3+p_4)}
\tilde e_{21}
\\
&+&
\biggl(
\frac{p_2 p_3}{
(p_1+p_2)(p_1+p_3)}
-
\frac{p_1p_3(p_1+p_2+p_3+p_4)}{(p_1+p_2)(p_1+p_3)(p_3+p_4)}
\biggr)
\tilde \varphi
\\
&+&
\biggl(
\frac{p_1 p_4}{
(p_1+p_3)(p_3+p_4)}
-
\frac{p_1p_3(p_1+p_2+p_3+p_4)}{(p_1+p_2)(p_1+p_3)(p_3+p_4)}
\biggr)
\tilde \phi,
\end{eqnarray*}
\begin{eqnarray*} 
 \phi &=&
\frac{p_4(p_2p_3-p_1p_4)}{
(p_1+p_2)(p_2+p_4)(p_3+p_4)}
\tilde e_{01}
\quad 
+ \quad 
\frac{p_2(p_1p_4-p_2p_3)}{
(p_1+p_2)(p_2+p_4)(p_3+p_4)}
\tilde e_{02}
\\
&+&
\frac{p_1p_2p_4(p_1+p_2+p_3+p_4)}{(p_1+p_2)(p_2+p_4)(p_2p_3-p_1p_4)}
\tilde e_{10}
\quad
- \quad
\frac{p_1 p_2}{
(p_1+p_2)(p_2+p_4)}
\tilde e_{12}
\\
&+&
\frac{p_2p_3p_4(p_1+p_2+p_3+p_4)}{(p_2+p_4)(p_3+p_4)(p_1p_4-p_2p_3)}
\tilde e_{20}
\quad - \quad
\frac{p_3 p_4}{
(p_2+p_4)(p_3+p_4)}
\tilde e_{21}
\\
&+&
\biggl(
\frac{p_1 p_4}{
(p_1+p_2)(p_2+p_4)}
-
\frac{p_2p_4(p_1+p_2+p_3+p_4)}{(p_1+p_2)(p_2+p_4)(p_3+p_4)}
\biggr)
\tilde \varphi
\\
&+&
\biggl(
\frac{p_2 p_3}{
(p_2+p_4)(p_3+p_4)}
-
\frac{p_2p_4(p_1+p_2+p_3+p_4)}{(p_1+p_2)(p_2+p_4)(p_3+p_4)}
\biggr) \tilde \phi.
\end{eqnarray*}
By an {\it antiautomorphism} of 
$\mathfrak{sl}_3(\C)$ we mean
a $\C$-linear bijection 
$\theta :
\mathfrak{sl}_3(\C) \to
\mathfrak{sl}_3(\C)$
such that $
\lbrack \beta,\gamma\rbrack^\theta = 
-\lbrack \beta^\theta,\gamma^\theta \rbrack
$
for all $\beta, \gamma \in 
\mathfrak{sl}_3(\C)$.
There exists an antiautomorphism $\dagger$ of 
$\mathfrak{sl}_3(\C)$ such that
$\beta^\dagger =\tilde W \beta^t \tilde W^{-1}
$ for all $\beta \in 
\mathfrak{sl}_3(\C)$.
Note that $(\beta^\dagger)^\dagger=\beta$ for all 
 $\beta \in 
\mathfrak{sl}_3(\C)$.
We have
\begin{eqnarray}
\rm
\begin{tabular}{c||ccc|ccc|cc}
$\beta$  &  $e_{01}$ & $e_{12}$ & $e_{02}$ & $e_{10}$ & $e_{21}$ & $e_{20}$
& $ \varphi $ & $ \phi  $
\\
\hline
$\beta^\dagger$ & 
$e_{10}\tilde \eta_1 /\tilde \eta_0$ 
&
$e_{21}\tilde \eta_2 /\tilde \eta_1$ 
& 
$e_{20}\tilde \eta_2 /\tilde \eta_0$ 
&  
$e_{01}\tilde \eta_0 /\tilde \eta_1$ 
&
$e_{12}\tilde \eta_1 /\tilde \eta_2$ 
&
$e_{02}\tilde \eta_0 /\tilde \eta_2$ 
&  
$\varphi$
&
$\phi $ 
\end{tabular}
\label{eq:table1}
\end{eqnarray}
Using $RWR^t = \tilde \vartheta \vartheta^{-1} \tilde W$
one finds
\begin{eqnarray}
\rm
\begin{tabular}{c||ccc|ccc|cc}
$\beta$  &  $\tilde e_{01}$ & $\tilde e_{12}$ & $\tilde e_{02}$ & $\tilde e_{10}$
& $\tilde e_{21}$ & $\tilde e_{20}$
& $\tilde \varphi$ & $\tilde \phi$
\\
\hline
$\beta^\dagger$ & 
$\tilde e_{10} \eta_1 /\eta_0$ 
&
$\tilde e_{21}\eta_2 / \eta_1$ 
& 
$\tilde e_{20}\eta_2 / \eta_0$ 
&  
$\tilde e_{01} \eta_0 / \eta_1$ 
&
$\tilde e_{12} \eta_1 / \eta_2$ 
&
$\tilde e_{02} \eta_0 / \eta_2$ 
&  
$\tilde \varphi $
&
$\tilde \phi $ 
\end{tabular}
\label{eq:table2}
\end{eqnarray}
Note that  $\dagger$ fixes each element
of $H$ and each element of
$\tilde H$. We now show that
$H$, $\tilde H$ together generate
$\mathfrak{sl}_3(\C)$. Define $\psi = -\varphi -\phi$
and 
$\tilde \psi = -\tilde \varphi -\tilde \phi$. Then
\begin{eqnarray*}
&&
e_{01} = \frac{-\lbrack \varphi,\lbrack \psi,\tilde \psi\rbrack\rbrack
+ \lbrack \varphi,\lbrack \varphi,\lbrack \psi,\tilde \psi\rbrack\rbrack
\rbrack}{2 \tilde \eta_0},
\qquad 
e_{10} = \frac{-\lbrack \varphi,\lbrack \psi,\tilde \psi\rbrack\rbrack
- \lbrack \varphi,\lbrack \varphi,\lbrack \psi,\tilde \psi\rbrack\rbrack
\rbrack}{2 \tilde \eta_1},
\\
&&
e_{02} = \frac{-\lbrack \phi,\lbrack \psi,\tilde \psi\rbrack\rbrack
+ \lbrack \phi,\lbrack \phi,\lbrack \psi,\tilde \psi\rbrack\rbrack
\rbrack}{2 \tilde \eta_0},
\qquad 
e_{20} = \frac{-\lbrack \phi,\lbrack \psi,\tilde \psi\rbrack\rbrack
- \lbrack \phi,\lbrack \phi,\lbrack \psi,\tilde \psi\rbrack\rbrack
\rbrack}{2 \tilde \eta_2},
\\
&&
e_{12} = \frac{-\lbrack \varphi,\lbrack \phi,\tilde \psi\rbrack\rbrack
- \lbrack \varphi,\lbrack \varphi,\lbrack \phi,\tilde \psi\rbrack\rbrack
\rbrack}{2 \tilde \eta_1},
\qquad 
e_{21} = \frac{-\lbrack \varphi,\lbrack \phi,\tilde \psi\rbrack\rbrack
+ \lbrack \varphi,\lbrack \varphi,\lbrack \phi,\tilde \psi\rbrack\rbrack
\rbrack}{2 \tilde \eta_2}.
\end{eqnarray*}

\section{An $\mathfrak{sl}_3(\C)$-module}\label{se3}
\noindent
We now define a certain
$\mathfrak{sl}_3(\C)$-module.
Let $x,y,z$ denote mutually commuting indeterminates.
Let $\C \lbrack x,y,z\rbrack$ denote the $\C$-algebra consisting
of the polynomials in $x,y,z$ that have all coefficients in
$\C$.
We abbreviate $A=
\C \lbrack x,y,z\rbrack$.
By a {\it derivation} of 
$A$
we mean a $\C$-linear map
$\partial: A \to A$ such that $\partial(\xi \zeta)=\partial(\xi) \zeta + 
\xi \partial (\zeta)$ for all $\xi,\zeta \in A$.
By 
\cite[p.~207]{carter}
the  space $A$
is an 
$\mathfrak{sl}_3(\C)$-module
on which each element of
$\mathfrak{sl}_3(\C)$
acts as a  derivation
and

\begin{center}
\begin{tabular}{c||ccc|ccc|cc}
$\xi$  &  $e_{01}.\xi$ & $e_{12}.\xi$ & $e_{02}.\xi$ & $e_{10}.\xi$
& $e_{21}.\xi$ & $e_{20}.\xi$
& $\varphi.\xi$ & $\phi.\xi$
\\
\hline
$x$ & 
$0$ & $0$ & $0$ &  
$y$ & $0$ & $z$ &  
$-x/3$ & $-x/3$ 
\\
$y$ &
$x$ & $0$ & $0$ &  
$0$ & $z$ & $0$ &  
$2y/3$ & $-y/3$ 
\\
$z$ &
$0$ & $y$ & $x$ &  
$0$ & $0$ & $0$ &  
$-z/3$ & $2z/3$ 
\end{tabular}
\end{center}
Let 
$V$ denote the subspace of $A$ consisting of the homogeneous polynomials
that have total degree $N$.
The following is a basis for $V$:
\begin{eqnarray}
x^r y^s z^t \qquad  r\geq 0, \quad s\geq 0, \quad t\geq 0, 
\qquad r+s+t = N.
\label{eq:Anbasis1-intro}
\end{eqnarray}
\noindent
The action of 
$\mathfrak{sl}_3(\C)$ on this basis is described as follows.

\begin{center}
\begin{tabular}{c||ccc}
$\xi$ &  $e_{01}.\xi$ & $e_{12}.\xi$ & $e_{02}.\xi$ 
\\
\hline
$x^r y^s z^t$ & 
$sx^{r+1}y^{s-1}z^t$ & $tx^r y^{s+1} z^{t-1}$ & $tx^{r+1} y^s z^{t-1}$   
\end{tabular}
\end{center}

\begin{center}
\begin{tabular}{c||ccc}
$\xi$ & $e_{10}.\xi$ & $e_{21}.\xi$ & $e_{20}.\xi$
\\
\hline
$x^r y^s z^t$ & 
$rx^{r-1} y^{s+1} z^t$ & $sx^r y^{s-1} z^{t+1}$ & $rx^{r-1} y^s z^{t+1}$
\end{tabular}
\end{center}

\begin{center}
\begin{tabular}{c||cc}
$\xi$ & $\varphi.\xi$ & $\phi.\xi$
\\
\hline
$x^r y^s z^t$ & 
$(s-N/3)x^r y^s z^t$ & $(t-N/3)x^r y^s z^t$ 
\end{tabular}
\end{center}
By the above data $V$ is an
$\mathfrak{sl}_3(\C)$-submodule of $A$, and this submodule
is irreducible
\cite[p.~97]{janzten}.
Let $\mathbb I$ denote
the set consisting of the 3-tuples of nonnegative integers
whose sum is $N$.
For $\lambda =(r,s,t) \in \mathbb I$
let $V_\lambda$ denote the subspace of
$V$ spanned by $x^ry^sz^t$. Then
$V = \sum_{\lambda \in \mathbb I} V_\lambda$
(direct sum)
and this is the weight space decomposition of $V$
with respect to $H$.
By construction
$\mbox{\rm dim}(V_\lambda)  = 1$ for all $\lambda \in \mathbb I$.

\medskip
\noindent We 
now consider the weight space decomposition
of $V$ with respect to $\tilde H$.
To describe this decomposition we make a change of variables.
Recall the matrix $R$ and
define
\begin{eqnarray*}
 \tilde x &=& R_{00}x+R_{10}y + R_{20}z,
\\
 \tilde y &=& R_{01}x+R_{11}y + R_{21}z,
\\ 
 \tilde z &=& R_{02}x+R_{12}y + R_{22}z.
\end{eqnarray*}
Thus $R$ is the transition matrix from
$x,y,z$ to $\tilde x,\tilde y,\tilde z$.
Using
$R=\tilde \vartheta \tilde W U^t$
we obtain
\begin{eqnarray}
\label{eq:keysum-intro}
\tilde \vartheta^{-1} \tilde x &=&
\tilde \eta_0 x +
\tilde \eta_1 y +
\tilde \eta_2 z,
\\
\tilde \vartheta^{-1} \tilde y  &=&
\tilde \vartheta^{-1} \tilde x - t \tilde \eta_1 y - v \tilde \eta_2 z,
\label{eq:handy1-intro}
\\
\tilde \vartheta^{-1} \tilde z  &=&
\tilde \vartheta^{-1} \tilde x - u \tilde \eta_1 y - w \tilde \eta_2 z.
\label{eq:handy2-intro}
\end{eqnarray}
Using $R^{-1}=\vartheta W U$ we obtain
\begin{eqnarray*}
 \vartheta^{-1} x &=&
 \eta_0 \tilde x +
 \eta_1 \tilde y +
 \eta_2 \tilde z,
\\
 \vartheta^{-1} y  &=&
 \vartheta^{-1}  x - t  \eta_1 \tilde y - u \eta_2 \tilde z,
\\
 \vartheta^{-1}  z  &=&
 \vartheta^{-1}  x - v  \eta_1 \tilde y - w \eta_2 \tilde z.
\end{eqnarray*}
The action of 
$\mathfrak{sl}_3(\C)$ on
$\tilde x$,
$\tilde y$,
$\tilde z$ is described as follows.

\begin{center}
\begin{tabular}{c||ccc|ccc|cc}
 $\xi$ &  $\tilde e_{01}.\xi$ & 
 $\tilde e_{12}.\xi$ & $ \tilde e_{02}.\xi$ & $\tilde e_{10}.\xi$ &
  $\tilde e_{21}.\xi$ & $\tilde e_{20}.\xi$
& $\tilde \varphi.\xi$ & $ \tilde \phi.\xi$
\\
\hline
$\tilde x$ & 
$0$ & $0$ & $0$ &  
$\tilde y$ & $0$ & $\tilde z$ &  
$-\tilde x/3$ & $-\tilde x/3$ 
\\
$\tilde y$ &
$\tilde x$ & $0$ & $0$ &  
$0$ & $\tilde z$ & $0$ &  
$2\tilde y/3$ & $-\tilde y/3$ 
\\
$\tilde z$ &
$0$ & $\tilde y$ & $\tilde x$ &  
$0$ & $0$ & $0$ &  
$-\tilde z/3$ & $2\tilde z/3$ 
\end{tabular}
\end{center}
The following is a basis for $V$.
\begin{eqnarray}
\tilde x^r \tilde y^s \tilde z^t \qquad  r\geq 0, \quad s\geq 0, 
\quad t\geq 0, \qquad r+s+t = N.
\label{eq:Anbasis2-intro}
\end{eqnarray}
The action of
$\mathfrak{sl}_3(\C)$ on 
this basis
 is described as follows.

\begin{center}
\begin{tabular}{c||ccc}
$\xi $ &  $\tilde e_{01}.\xi$ & $\tilde e_{12}.\xi$ & $\tilde e_{02}.\xi$ 
\\
\hline
$\tilde x^r \tilde y^s \tilde z^t$ & 
$s\tilde x^{r+1}\tilde y^{s-1} \tilde z^t$ & $t\tilde x^r \tilde y^{s+1}
\tilde z^{t-1}$ & $t \tilde x^{r+1} \tilde y^s \tilde z^{t-1}$   
\end{tabular}
\end{center}

\begin{center}
\begin{tabular}{c||ccc}
$\xi$ & $\tilde e_{10}.\xi$ & $\tilde e_{21}.\xi$ & $\tilde e_{20}.\xi $
\\
\hline
$\tilde x^r \tilde y^s \tilde  z^t$ & 
$r \tilde x^{r-1} \tilde y^{s+1} \tilde z^t$ & $s\tilde x^r \tilde y^{s-1} \tilde  z^{t+1}$ & $ r \tilde x^{r-1} \tilde  y^s  \tilde z^{t+1}$
\end{tabular}
\end{center}

\begin{center}
\begin{tabular}{c||cc}
$\xi$& $\tilde \varphi.\xi$ & $\tilde \phi.\xi$
\\
\hline
$\tilde x^r \tilde y^s \tilde z^t$ & 
$(s-N/3)\tilde x^r \tilde  y^s \tilde  z^t$ & 
$(t-N/3)\tilde x^r \tilde y^s \tilde z^t$ 
\end{tabular}
\end{center}

\noindent
For each $\lambda =(r,s,t)\in \mathbb I$ let
${\tilde V}_\lambda$ denote the subspace of
$V$ spanned by
$\tilde x^r \tilde y^s \tilde z^t$.
Observe that $V=\sum_{\lambda \in \mathbb I} {\tilde V}_\lambda$
(direct sum) and this is the weight space decomposition of
$V$ with respect to $\tilde H$.
By construction
$\mbox{\rm dim}(\tilde V_\lambda)  = 1$ for all
$\lambda \in \mathbb I$.

\medskip
\noindent 
We comment on how $H$ and $\tilde H$ act 
on the  weight spaces of the other one.
A pair of elements $(r,s,t)$ and $(r',s',t')$ in $\mathbb I$
will be called {\it adjacent} whenever
$(r-r',s-s',t-t')$ is a permutation of
$(1,-1,0)$.
Then $H$ and $\tilde H$ act on the weight spaces of the
other one
as follows.
For all $\lambda \in \mathbb I$,
\begin{eqnarray}
\label{eq:mutual}
\tilde H V_\lambda \subseteq  V_\lambda + 
\sum_{\genfrac{}{}{0pt}{} {\mu \in \mathbb I}{\mu \;{\rm adj}\;\lambda}}
V_\mu,
\qquad \qquad 
H \tilde V_\lambda \subseteq  \tilde V_\lambda + 
\sum_{\genfrac{}{}{0pt}{} {\mu \in \mathbb I}{\mu \;{\rm adj}\;\lambda}}
\tilde V_\mu.
\end{eqnarray}

\section{A bilinear form}\label{se4}
\noindent In this section we introduce a
symmetric bilinear form
$\langle \,,\,\rangle $ on
$V$. As we will see,
both
\begin{eqnarray}
\langle V_\lambda, V_\mu\rangle &=& 0
\quad {\rm if} \quad \lambda \not = \mu,
\qquad \qquad \lambda, \mu \in \mathbb I,
\label{eq:orth1}
\\
\langle \tilde V_\lambda, \tilde V_\mu\rangle &=& 0
\quad {\rm if} \quad \lambda \not = \mu,
\qquad \qquad \lambda, \mu \in \mathbb I.
\label{eq:orth2}
\end{eqnarray}
We define  
$\langle \,,\,\rangle$ as follows.
With respect to $\langle \,,\,\rangle$ 
the vectors
 (\ref{eq:Anbasis1-intro}) are mutually orthogonal
 and
\begin{eqnarray}
\Vert x^r y^s z^t \Vert^2 = \frac{r!s!t!}
{\tilde \eta^r_0
\tilde \eta^s_1
\tilde \eta^t_2}
\vartheta^N
\qquad \qquad
r\geq 0, 
\quad s\geq 0,
\quad t\geq 0,
\quad r+s+t=N.
\label{eq:basis1-norm-intro}
\end{eqnarray}
We are using the notation 
$\Vert \xi\Vert^2=\langle \xi,\xi \rangle$.
The form $\langle \,,\,\rangle $ is symmetric,
nondegenerate, and satisfies
(\ref{eq:orth1}).
Using
(\ref{eq:table1}) and the
tables below
(\ref{eq:Anbasis1-intro}) we obtain
\begin{eqnarray}
\langle \beta \xi,\zeta\rangle = \langle \xi,\beta^\dagger \zeta\rangle
\qquad \qquad 
\forall \beta \in
\mathfrak{sl}_3(\C), \qquad   \forall \xi,\zeta \in V.
\label{lem:bilanti-intro}
\end{eqnarray}
Line 
 (\ref{eq:orth2}) follows 
from 
(\ref{lem:bilanti-intro})
and since
$\dagger$ fixes each element of $\tilde H$.
The following is the basis for $V$ that is
 dual to
 (\ref{eq:Anbasis1-intro}) with respect to
$\langle \,,\,\rangle$.
\begin{eqnarray}
 \frac
 {
\tilde \eta^r_0
\tilde \eta^s_1
\tilde \eta^t_2}
 {r!s!t!}
\,\frac{x^r  y^s  z^t}{\vartheta^{N}} \qquad \qquad
r\geq 0, \quad s\geq 0, \quad t\geq 0, \qquad 
r+s+t = N.
\label{eq:dualbasis1-intro}
\end{eqnarray}
 The vector $\tilde x^N$ is equal to
 $N!\nu^N$ times the sum
of the 
vectors (\ref{eq:dualbasis1-intro});
this is verified using
(\ref{eq:keysum-intro}) and the trinomial theorem.

\begin{lemma}
\label{lem:dualnorm-intro}
With respect to $\langle \,,\,\rangle $
the vectors
 (\ref{eq:Anbasis2-intro}) are mutually orthogonal
and
\begin{eqnarray}
\Vert \tilde x^r \tilde y^s \tilde z^t \Vert^2 = \frac{r!s!t!}
{ \eta^r_0
 \eta^s_1
 \eta^t_2}
\tilde \vartheta^N
\qquad \qquad
r\geq 0, 
\quad s\geq 0,
\quad t\geq 0,
\quad r+s+t=N.
\label{eq:dualsquare-intro}
\end{eqnarray}
\end{lemma}
\noindent {\it Proof:} 
By 
 (\ref{eq:orth2})
the vectors
(\ref{eq:Anbasis2-intro})  are mutually
orthogonal.
We now verify
(\ref{eq:dualsquare-intro}).
We do this by induction
on $s+t$.
First assume $s+t=0$, so that $(r,s,t)=(N,0,0)$.
We must show
$\Vert \tilde x^N\Vert^2 = N!\, 
\eta_0^{-N}
\tilde \vartheta ^N $. 
The vector ${\tilde x}^N$ is equal to
$N!\nu^N$ times the sum of the vectors
(\ref{eq:dualbasis1-intro}), and these vectors
are mutually orthogonal.
Therefore  
$\Vert \tilde x^N\Vert^2$
is equal to
$(N!)^2 \nu^{2N}$ times the sum of the square norms of
the vectors 
(\ref{eq:dualbasis1-intro}).
Computing this sum using
(\ref{eq:basis1-norm-intro}) we obtain
\begin{eqnarray*}
\Vert \tilde x^N\Vert^2 
&=& 
(N!)^2\,\nu^{2N}\vartheta^{-N} 
\sum_{\genfrac{}{}{0pt}{} {0\leq r,s,t}{r+s+t=N}}
 \frac
 {
\tilde \eta^r_0
\tilde \eta^s_1
\tilde \eta^t_2}
 {r!s!t!}
\\
&=& 
N!\,\nu^{2N} \vartheta^{-N} 
\sum_{\genfrac{}{}{0pt}{} {0\leq r,s,t}{r+s+t=N}}
\tilde \eta^r_0
\tilde \eta^s_1
\tilde \eta^t_2
\binom{N}{r\;s\;t}
\\
&=& 
N!\,\nu^{2N}\vartheta^{-N} 
(
\tilde \eta^r_0 +
\tilde \eta^s_1 +
\tilde \eta^t_2)^N
\\
&=& 
N!\,\nu^{2N}\vartheta^{-N} 
\\
&=& 
N!\,
\eta_0^{-N}
\tilde \vartheta^{N}.
\end{eqnarray*}
%We have verified 
%(\ref{eq:dualsquare-intro})
%for $s+t=0$.
Next assume
$s+t>0$, so that $s>0$ or
$t>0$. We assume $s>0$; the case
$t>0$ is similar.
By
(\ref{lem:bilanti-intro}), 
\begin{eqnarray}
\label{eq:ee}
\langle
\tilde e_{01}.\tilde x^r \tilde y^s \tilde z^t ,\tilde x^{r+1} \tilde y^{s-1} \tilde z^t\rangle
=
\langle
\tilde x^r \tilde y^s \tilde z^t , \tilde e^\dagger_{01}. \tilde x^{r+1} \tilde y^{s-1} 
\tilde z^t\rangle.
\end{eqnarray}
Using the first table below
(\ref{eq:Anbasis2-intro}) and then induction,
\begin{eqnarray*}
\langle
\tilde e_{01}.\tilde x^r \tilde y^s \tilde z^t ,\tilde x^{r+1} \tilde y^{s-1} \tilde z^t\rangle
&=&
\Vert 
\tilde x^{r+1} \tilde y^{s-1} \tilde z^t
\Vert^2 s
\\
&=&
\frac{(r+1)!s!t!}
{ \eta^{r+1}_0
 \eta^{s-1}_1
 \eta^t_2}
\tilde \vartheta^N.
\end{eqnarray*}
Using
(\ref{eq:table2}) 
and the second table below
(\ref{eq:Anbasis2-intro}),
\begin{eqnarray*}
\langle
\tilde x^r \tilde y^s \tilde z^t , \tilde e^\dagger_{01} .\tilde x^{r+1}
\tilde y^{s-1} 
\tilde z^t\rangle
&=&
\langle
\tilde x^r \tilde y^s \tilde z^t , \tilde e_{10} .\tilde x^{r+1}
\tilde y^{s-1} 
\tilde z^t\rangle \eta_1/\eta_0
\\
&=&
\Vert 
\tilde x^r \tilde y^s \tilde z^t\Vert^2 (r+1)\eta_1/\eta_0.
\end{eqnarray*}
Evaluating
(\ref{eq:ee}) using these comments we obtain
\begin{eqnarray*}
\Vert \tilde x^r \tilde y^s \tilde z^t \Vert^2 = \frac{r!s!t!}
{ \eta^r_0
 \eta^s_1
 \eta^t_2}
\tilde \vartheta^N.
\end{eqnarray*}
We have verified
(\ref{eq:dualsquare-intro}) and the lemma is proved.
\hfill $\Box$ \\

\noindent
The following is the basis for $V$ that is
 dual to
 (\ref{eq:Anbasis2-intro}) with respect to
$\langle \,,\,\rangle$.
\begin{eqnarray}
 \frac
 {
 \eta^r_0
 \eta^s_1
 \eta^t_2}
 {r!s!t!}
\,\frac{\tilde x^r \tilde y^s  \tilde z^t}{\tilde \vartheta^{N}} \qquad \qquad  
r\geq 0, \quad s\geq 0, \quad t\geq 0, \qquad
r+s+t=N.
\label{eq:dualbasis2-intro}
\end{eqnarray}
The vector
 $\tilde x^N$ is equal to
 $N!\nu^N$ times the sum
of the 
vectors (\ref{eq:dualbasis2-intro}).

\section{The Rahman polynomials and  $\mathfrak{sl}_3(\C)$}\label{se5} 
\medskip
\noindent
%We now come to our main results, in which 
%the Rahman polynomials are related to the
%$\mathfrak{sl}_3(\C)$-module
%$V$.
In this section 
the Rahman polynomials are related to the
$\mathfrak{sl}_3(\C)$-module
$V$.
%%c6-14

\medskip
\noindent
We would like to acknowledge that the 
following result is similar to
\cite[Section~2]{tanaka}, although the setup
is somewhat different.

\begin{theorem}
\label{thm:main1}
For a basis vector 
$\tilde x^\rho \tilde y^\sigma \tilde z^\tau$ from
(\ref{eq:Anbasis2-intro}),
\begin{eqnarray}
\tilde x^\rho \tilde y^\sigma \tilde z^\tau
&=&
N! \nu^N  
%\tilde \lambda^N
\sum_{\genfrac{}{}{0pt}{} {0\leq r,s,t}{r+s+t=N}}
P (s,t,\sigma,\tau)
%\binom{N}{r \; \;s\;\; t} 
\,\frac{\tilde\eta^r_0
\tilde \eta^s_1
\tilde \eta^t_2}{r!s!t!}
\,\frac{
x^r y^s z^t}{\vartheta^N}.
\label{eq:trans2-intro}
\end{eqnarray}
For a basis vector
$x^\rho y^\sigma z^\tau$ from
(\ref{eq:Anbasis1-intro}),
\begin{eqnarray}
x^\rho  y^\sigma  z^\tau
&=&
N! \nu^N
%\lambda^N
\sum_{\genfrac{}{}{0pt}{} {0\leq r,s,t}{r+s+t=N}}
 P (\sigma,\tau,s,t)
%\binom{N}{r \;\; s\;\; t} 
\,\frac{\eta^r_0
 \eta^s_1
 \eta^t_2}{r!s!t!}
\,\frac{\tilde x^r \tilde y^s\tilde z^t}{\tilde \vartheta^N}.
\label{eq:trans1-intro}
\end{eqnarray}
\end{theorem}
\noindent {\it Proof:} 
We first prove 
(\ref{eq:trans2-intro}).
Since $\rho+\sigma+\tau=N$,
\begin{eqnarray}
\label{eq:view-intro}
\tilde x^\rho \tilde y^\sigma \tilde z^\tau
= {\tilde \vartheta}^N 
(\tilde \vartheta^{-1} \tilde x)^\rho
(\tilde \vartheta^{-1} \tilde y)^\sigma
(\tilde \vartheta^{-1} \tilde z)^\tau.
\end{eqnarray}
Using
(\ref{eq:handy1-intro}) and the trinomial theorem,
\begin{eqnarray*}
(\tilde \vartheta^{-1} \tilde y)^\sigma &=&
(\tilde \vartheta^{-1} \tilde x - t \tilde \eta_1 y - 
v \tilde \eta_2 z)^\sigma
\\
&=&
\sum_{\genfrac{}{}{0pt}{} {0\leq i,k}{i+k\leq \sigma}}
(\tilde \vartheta^{-1} \tilde x)^{\sigma-i-k}
(-t\tilde \eta_1 y)^i
(-v\tilde \eta_2 z)^k
\binom{\sigma}{i \;\; k \;\; \sigma-i-k} 
\\
&=&
\sum_{\genfrac{}{}{0pt}{} {0\leq i,k}{i+k\leq \sigma}}
(\tilde \vartheta^{-1} \tilde x)^{\sigma-i-k}
(t\tilde \eta_1 y)^i
(v\tilde \eta_2 z)^k
\frac{(-\sigma)_{i+k}}{i!k!}.
\end{eqnarray*}
Similarly using
(\ref{eq:handy2-intro}),
\begin{eqnarray*}
(\tilde \vartheta^{-1} \tilde z)^\tau &=&
(\tilde \vartheta^{-1} \tilde x - u \tilde \eta_1 y - 
w \tilde \eta_2 z)^\tau
\\
&=&
\sum_{\genfrac{}{}{0pt}{} {0\leq j,\ell}{j+\ell\leq \tau}}
(\tilde \vartheta^{-1} \tilde x)^{\tau-j-\ell}
(-u\tilde \eta_1 y)^j
(-w\tilde \eta_2 z)^\ell
\binom{\tau}{j \;\; \ell \;\; \tau-j-\ell} 
\\
&=&
\sum_{\genfrac{}{}{0pt}{} {0\leq j,\ell}{j+\ell\leq \tau}}
(\tilde \vartheta^{-1} \tilde x)^{\tau-j-\ell}
(u\tilde \eta_1 y)^j
(w\tilde \eta_2 z)^\ell
\frac{(-\tau)_{j+\ell}}{j!\ell !}.
\end{eqnarray*}
Evaluating
(\ref{eq:view-intro}) using the above comments
and 
(\ref{eq:principle-intro}), we 
 obtain
\begin{eqnarray}
\tilde x^\rho \tilde y^\sigma \tilde z^\tau
&=& 
{\tilde \vartheta}^N 
\sum_{\genfrac{}{}{0pt}{} {0\leq i,j,k,\ell}{i+j+k+\ell \leq N}}
t^i u^j v^k w^\ell \frac{(-\sigma)_{i+k} (-\tau)_{j+\ell}}
{i!j!k!\ell!} \quad \times 
\nonumber
\\
&&
\qquad (\tilde \vartheta^{-1} \tilde x)^{N-i-j-k-\ell}
(\tilde \eta_1 y)^{i+j}
(\tilde \eta_2 z)^{k+\ell}.
\label{eq:storm-intro}
\end{eqnarray}
In the expression
(\ref{eq:storm-intro}) consider
the first factor. 
Using
(\ref{eq:keysum-intro}) and the trinomial theorem,
\begin{eqnarray*}
(\tilde \vartheta^{-1} \tilde x)^{N-i-j-k-\ell}
&=& (\tilde \eta_0 x +
\tilde \eta_1 y +
\tilde \eta_2 z)^{N-i-j-k-\ell}
\\
&=&
 \sum_{\genfrac{}{}{0pt}{} {0\leq r, b,c}{r+b+c=N-i-j-k-\ell}}
(\tilde \eta_0 x)^r 
(\tilde \eta_1 y)^{b}
(\tilde \eta_2 z)^{c}
\binom{N-i-j-k-\ell}{r\; \;\; b \;\;\; c}.
\end{eqnarray*}
Therefore 
(\ref{eq:storm-intro}) 
is equal to
\begin{eqnarray*}
 \sum_{\genfrac{}{}{0pt}{} {0\leq r, b,c}{r+b+c=N-i-j-k-\ell}}
(\tilde \eta_0 x)^r 
(\tilde \eta_1 y)^{b+i+j}
(\tilde \eta_2 z)^{c+k+\ell}
\binom{N-i-j-k-\ell}{r \;\; \;b \;\;\; c},
\end{eqnarray*}
which is equal to
\begin{eqnarray*} 
&& \sum_{\genfrac{}{}{0pt}{} {0\leq r, b,c}{r+b+c=N-i-j-k-\ell}}
(\tilde \eta_0 x)^r 
(\tilde \eta_1 y)^{b+i+j}
(\tilde \eta_2 z)^{c+k+\ell}
\quad \times
\\
&& 
\qquad \qquad
\qquad \qquad
\frac{(-b-i-j)_{i+j} (-c-k-\ell)_{k+\ell}}{(-N)_{i+j+k+\ell}}
\binom{N}{r \;\; \;b+i+j \;\;\; c+k+\ell}.
\end{eqnarray*}
After changing variables
$s=b+i+j$,
$t=c+k+\ell$ and using
(\ref{eq:principle-intro})
the above sum becomes
\begin{eqnarray}
\label{eq:sumfin-intro}
\sum_{\genfrac{}{}{0pt}{} {0\leq r, s,t}{r+s+t=N}}
(\tilde \eta_0 x)^r 
(\tilde \eta_1 y)^s
(\tilde \eta_2 z)^t
\frac{(-s)_{i+j} (-t)_{k+\ell}}{(-N)_{i+j+k+\ell}}
\binom{N}{r\;\;s\;\;t}.
\end{eqnarray}
Therefore  (\ref{eq:storm-intro}) is equal to
(\ref{eq:sumfin-intro}).
Upon replacing 
(\ref{eq:storm-intro}) by
(\ref{eq:sumfin-intro})
the expression
$\tilde x^\rho \tilde y^\sigma \tilde z^\tau$
becomes a sum which is equal to
the right-hand side of
(\ref{eq:trans2-intro}).
This proves
(\ref{eq:trans2-intro})
and the proof of
(\ref{eq:trans1-intro}) is similar.
\hfill $\Box$ \\

\begin{theorem}
For nonnegative integers $s,t$ whose sum is at most $N$, both
\begin{eqnarray}
\label{eq:p1-intro}
 P (s,t, \tilde \varphi+N/3 \, I, \tilde \phi + N/3 \,I) x^N
&=& x^r y^s z^t,
\\
\label{eq:p2-intro}
P (\varphi+N/3 \,I, \phi + N/3 \,I ,s,t) \tilde x^N
&=& \tilde x^r \tilde y^s \tilde z^t,
\end{eqnarray}
where $r=N-s-t$.
\end{theorem}
\noindent {\it Proof:} 
Concerning 
(\ref{eq:p2-intro}),
 observe
\begin{eqnarray*}
\tilde x^r \tilde y^s \tilde z^t
&=& 
N! \nu^N
%\tilde \lambda^N 
\sum_{\genfrac{}{}{0pt}{} {0\leq i,j,k}{i+j+k=N}}
 P (j,k,s,t)
%\binom{N}{i \; j\; k}
\,\frac{
\tilde \eta^i_0 
\tilde \eta^j_1 
\tilde \eta^k_2 }{i!j!k!}
\,
\frac{
x^i y^j z^k}{\vartheta^N}
\\
&=& 
 P (\varphi+N/3\,I,\phi+N/3\,I,s,t)
%\tilde \lambda^N 
N! \nu^N
\sum_{\genfrac{}{}{0pt}{} {0\leq i,j,k}{i+j+k=N}}
%\binom{N}{i \; j\; k}
\,\frac{\tilde \eta^i_0 
\tilde \eta^j_1 
\tilde \eta^k_2 }{i!j!k!}
\,
\frac{
x^i y^j z^k}{\vartheta^N}
\\
&=&
 P (\varphi+N/3\,I,\phi+N/3\,I,s,t)\tilde x^N.
\end{eqnarray*}
Line
(\ref{eq:p1-intro}) is similarly
obtained.
\hfill $\Box$ \\

\begin{theorem}
\label{thm:pcosines-intro}
For a basis vector $x^ry^sz^t$ from
(\ref{eq:Anbasis1-intro}) and
a basis vector
$\tilde x^\rho \tilde y^\sigma \tilde z^\tau$ from
(\ref{eq:Anbasis2-intro}),
\begin{eqnarray*}
\langle 
x^ry^sz^t,
\tilde x^\rho \tilde y^\sigma \tilde z^\tau
\rangle = N! \, \nu^N   P (s,t,\sigma,\tau).
\end{eqnarray*}
\end{theorem}
\noindent {\it Proof:} 
By 
(\ref{eq:trans2-intro}),
\begin{eqnarray*}
\tilde x^\rho \tilde y^\sigma \tilde z^\tau
&=&
%\tilde \lambda^N
N! \nu^N
\sum_{\genfrac{}{}{0pt}{} {0\leq i,j,k}{i+j+k=N}}
P (j,k,\sigma,\tau)
%\binom{N}{i \; \;j\;\; k} 
\,
\frac{
\tilde\eta^i_0
\tilde \eta^j_1
\tilde \eta^k_2}{i!j!k!}
\,\frac{
x^i y^j z^k}{\vartheta^N}.
\end{eqnarray*}
In the above equation take the
inner product of
each side with 
$x^r y^s z^t$,
and use the fact that
the bases 
(\ref{eq:Anbasis1-intro}),
(\ref{eq:dualbasis1-intro})
are dual
with respect to
$\langle \,,\,\rangle$.
\hfill $\Box$ \\

\noindent
\noindent {\it Proof of Theorem
\ref{thm:rahmanorthog}}.
We first verify
(\ref{eq:orthog1-intro}).
Abbreviate $\rho= N-\sigma-\tau$.
By
(\ref{eq:trans2-intro}) we have
\begin{eqnarray*}
\tilde x^\rho \tilde y^\sigma \tilde z^\tau
&=&
%\tilde \lambda^N
N! \nu^N
\sum_{\genfrac{}{}{0pt}{} {0\leq i,j,k}{i+j+k=N}}
 P (j,k,\sigma,\tau)
%\binom{N}{i \; \;j\;\; k} 
\,
\frac{
\tilde\eta^i_0
\tilde \eta^j_1
\tilde \eta^k_2}{i!j!k!}
\,
\frac{x^i y^j z^k}{\vartheta^N}.
\end{eqnarray*}
In the above equation take the
inner product of
each side with 
$\tilde x^r \tilde y^s  \tilde z^t$.
In the resulting equation
simplify the left-hand side
using
Lemma
\ref{lem:dualnorm-intro}
and simplify
the right-hand side using
Theorem
\ref{thm:pcosines-intro}.
This gives
(\ref{eq:orthog1-intro}).
Line
(\ref{eq:orthog2-intro}) is similarly verified.
\hfill $\Box$ \\

\section{Some seven-term recurrence relations}

In this section we
display
two seven-term recurrence relations 
satisfied by the 
Rahman polynomials,
along with the corresponding relations satisfied
by their duals.
 To obtain these relations
we use the
$\mathfrak{sl}_3(\C)$-module $V$ from Section 3.
We would like to acknowledge here the
work of
Gr\"unbaum 
\cite[Section~6]{grun}, which hints at
the existence of these relations
and contains 
a wealth of related data 
for the case $N=5$.

\begin{theorem}
\label{thm:7termrec}
Fix nonnegative integers $s$, $t$ whose sum is at most $N$,
and nonnegative integers $\sigma$, $ \tau$ whose sum is at most $N$.
Then {\rm (i)--(iv)} hold below.
\begin{enumerate}
\item[\rm (i)]
$(s-N/3)P(s,t,\sigma,\tau)$ is a weighted sum
with the following terms and coefficients:

\begin{tabular}{c|c}
{\rm term}  &  {\rm coefficient} 
\\
\hline
%\end{tabular}
$P(s,t,\sigma-1,\tau)$
&
$
\sigma \, \frac{p_3(p_1p_4-p_2p_3)}{
(p_1+p_2)(p_1+p_3)(p_3+p_4)}
$
%\tilde e_{01}
%\sigma P(s,t,\sigma-1,\tau)
\\
%\quad 
%+ \quad
$ 
P(s,t,\sigma,\tau-1)
$
&
$
\tau \,\frac{p_1(p_2p_3-p_1p_4)}{
(p_1+p_2)(p_1+p_3)(p_3+p_4)}
$
%\tilde e_{02}
%\tau P(s,t,\sigma,\tau-1)
\\
$
P(s,t,\sigma+1,\tau)
$
&
$
\rho
\,
\frac{p_1p_2p_3(p_1+p_2+p_3+p_4)}{(p_1+p_2)(p_1+p_3)(p_1p_4-p_2p_3)}
$
%\tilde e_{10}
%\rho P(s,t,\sigma+1,\tau)
%\quad -
%\quad
\\
$
P(s,t,\sigma+1,\tau-1)
$
&
$
\tau
\,\frac{-p_1 p_2}{
(p_1+p_2)(p_1+p_3)}
$
%\tilde e_{12}
%\tau P(s,t,\sigma+1,\tau-1)
\\
$
P(s,t,\sigma,\tau+1)
$
& 
$
\rho \,
\frac{p_1p_3p_4(p_1+p_2+p_3+p_4)}{(p_1+p_3)(p_3+p_4)(p_2p_3-p_1p_4)}
$
%\tilde e_{20}
%\rho P(s,t,\sigma,\tau+1)
%\quad - \quad
\\
$
P(s,t,\sigma-1,\tau+1)
$
&
$
\sigma \,\frac{-p_3 p_4}{
(p_1+p_3)(p_3+p_4)}
$
%\tilde e_{21}
%\sigma P(s,t,\sigma-1,\tau+1)
\\
%&& + \quad
$
P(s,t,\sigma,\tau)
$
&
$
(\sigma-N/3)
%\biggl(
\Bigl(
\frac{p_2 p_3}{
(p_1+p_2)(p_1+p_3)}
-
\,
\frac{p_1p_3(p_1+p_2+p_3+p_4)}{(p_1+p_2)(p_1+p_3)(p_3+p_4)}
\Bigr)
%\tilde \varphi
$
\\
&
$
+\; (\tau-N/3)
\Bigl(
\frac{p_1 p_4}{
(p_1+p_3)(p_3+p_4)}
-
\frac{p_1p_3(p_1+p_2+p_3+p_4)}{(p_1+p_2)(p_1+p_3)(p_3+p_4)}
\Bigr)
$
%P(s,t,\sigma,\tau).
%\tilde \phi,
\end{tabular}
\item[\rm (ii)]
$(t-N/3)P(s,t,\sigma,\tau)$ is
a weighted sum with the following terms and coefficients:

\begin{tabular}{c|c}
{\rm term}  &  {\rm coefficient} 
\\
\hline
$
P(s,t,\sigma-1,\tau)
$
&
$
\sigma \,\frac{p_4(p_2p_3-p_1p_4)}{
(p_1+p_2)(p_2+p_4)(p_3+p_4)}
%\sigma E_{\sigma}^{-1}
$
\\
$
P(s,t,\sigma,\tau-1)
$
&
$
\tau \,
\frac{p_2(p_1p_4-p_2p_3)}{
(p_1+p_2)(p_2+p_4)(p_3+p_4)}
%\tau E_{\tau}^{-1}
$
\\
$
P(s,t,\sigma+1,\tau)
$
&
$
\rho \,\frac{p_1p_2p_4(p_1+p_2+p_3+p_4)}{(p_1+p_2)(p_2+p_4)(p_2p_3-p_1p_4)}
%(N-\sigma-\tau)E_{\sigma}
%- \quad
$
\\
$
P(s,t,\sigma+1,\tau-1)
$
&
$
\tau
\,
\frac{-p_1 p_2}{
(p_1+p_2)(p_2+p_4)}
%\tau E_{\sigma}E_{\tau}^{-1}
$
\\
$
P(s,t,\sigma,\tau+1)
$
&
$
\rho
\,
\frac{p_2p_3p_4(p_1+p_2+p_3+p_4)}{(p_2+p_4)(p_3+p_4)(p_1p_4-p_2p_3)}
%(N-\sigma-\tau)E_{\tau}
%\quad - \quad
$
\\
$
P(s,t,\sigma-1,\tau+1)
$
&
$
\sigma
\,
\frac{-p_3 p_4}{
(p_2+p_4)(p_3+p_4)}
%\sigma E_{\sigma}^{-1}E_{\tau}
$
\\
%&+&
$
P(s,t,\sigma,\tau)
$
&
$
(\sigma-N/3)
\Bigl(
\frac{p_1 p_4}{
(p_1+p_2)(p_2+p_4)}
-
\frac{p_2p_4(p_1+p_2+p_3+p_4)}{(p_1+p_2)(p_2+p_4)(p_3+p_4)}
\Bigr)
%(\sigma-N/3)\,\rm{Id}
$
\\
&$+\;
(\tau-N/3)
\Bigl(
\frac{p_2 p_3}{
(p_2+p_4)(p_3+p_4)}
-
\frac{p_2p_4(p_1+p_2+p_3+p_4)}{(p_1+p_2)(p_2+p_4)(p_3+p_4)}
\Bigr)$
\end{tabular}
\newpage
\item[\rm (iii)]
$(\sigma - N/3)P(s,t,\sigma,\tau)$
is a weighted sum with the following terms and coefficients:

\begin{tabular}{c|c}
{\rm term}  &  {\rm coefficient} 
\\
\hline
$
P(s-1,t,\sigma,\tau)
$
&
$
s \,\frac{p_2(p_1p_4-p_2p_3)}{
(p_1+p_2)(p_1+p_3)(p_2+p_4)}
$
%sE_s^{-1}
%\quad
%+
%\quad
\\
$
P(s,t-1,\sigma,\tau)
$
&
$
t\,
\frac{p_1(p_2p_3-p_1p_4)}{
(p_1+p_2)(p_1+p_3)(p_2+p_4)} 
$
%tE_t^{-1}
\\
$
P(s+1,t,\sigma,\tau)
$
&
$
r \,
\frac{p_1p_2p_3(p_1+p_2+p_3+p_4)}{(p_1+p_2)(p_1+p_3)(p_1p_4-p_2p_3)}
$
%(N-s-t)E_s
%\quad
%-
%\quad
\\
$
P(s+1,t-1,\sigma,\tau)
$
&
$
t 
\,
\frac{-p_1 p_3}{(p_1+p_2)(p_1+p_3)}
$
%tE_sE_t^{-1}
\\
$
P(s,t+1,\sigma,\tau)
$
&
$
r 
\,
\frac{p_1p_2p_4(p_1+p_2+p_3+p_4)}{(p_1+p_2)(p_2+p_4)(p_2p_3-p_1p_4)}
$
%(N-s-t)E_t
%\quad
%-
%\quad
\\
$
P(s-1,t+1,\sigma,\tau)
$
&
$
s
\,
\frac{-p_2 p_4}{
(p_1+p_2)(p_2+p_4)}
%sE_s^{-1}E_t
$
\\
$
P(s,t,\sigma,\tau)
$
&
$
(s-N/3)
\Bigl(\frac{p_2 p_3}{
(p_1+p_2)(p_1+p_3)}
-
\frac{p_1p_2(p_1+p_2+p_3+p_4)}{(p_1+p_2)(p_1+p_3)(p_2+p_4)}
\Bigr) 
$
\\
& 
$+\;
(t-N/3)
\Bigl(\frac{p_1 p_4}{
(p_1+p_2)(p_2+p_4)}
-
\frac{p_1p_2(p_1+p_2+p_3+p_4)}{(p_1+p_2)(p_1+p_3)(p_2+p_4)}
\Bigr)
$
\end{tabular}

\item[\rm (iv)]
$(\tau-N/3)P(s,t,\sigma,\tau)$ is a weighted sum with the
following terms and coefficients:

\begin{tabular}{c|c}
{\rm term}  &  {\rm coefficient} 
\\
\hline
$
P(s-1,t,\sigma,\tau)
$
&
$
s
\,
\frac{p_4(p_2p_3-p_1p_4)}{
(p_1+p_3)(p_3+p_4)(p_2+p_4)}
%sE_s^{-1}
$
\\
$
P(s,t-1,\sigma,\tau)
$
&
$
t
\,
\frac{p_3(p_1p_4-p_2p_3)}{
(p_1+p_3)(p_3+p_4)(p_2+p_4)}$
%tE_t^{-1}
\\
$
P(s+1,t,\sigma,\tau)
$
&
$
r
\,
\frac{p_1p_3p_4(p_1+p_2+p_3+p_4)}{(p_1+p_3)(p_3+p_4)(p_2p_3-p_1p_4)}$
%(N-s-t)E_s
%\quad
\\
$
P(s+1,t-1,\sigma,\tau)
$
&
$
t
\,
%- 
%\quad
\frac{-p_1 p_3}{
(p_1+p_3)(p_3+p_4)}%
%tE_sE_t^{-1}
$
\\
$
P(s,t+1,\sigma,\tau)
$
&
$
r
\,
\frac{p_2p_3p_4(p_1+p_2+p_3+p_4)}{(p_2+p_4)(p_3+p_4)(p_1p_4-p_2p_3)}$
%(N-s-t)E_t
%\quad
%-
%\quad
\\
$
P(s-1,t+1,\sigma,\tau)
$
&
$
s
\,
\frac{-p_2 p_4}{
(p_2+p_4)(p_3+p_4)}$
%sE_s^{-1}E_t
\\
$
P(s,t,\sigma,\tau)
$
&
$
(s-N/3)
\Bigl(
\frac{p_1 p_4}{
(p_1+p_3)(p_3+p_4)}
-
\frac{p_3p_4(p_1+p_2+p_3+p_4)}{(p_1+p_3)(p_3+p_4)(p_2+p_4)}
\Bigr)
%(s-N/3)\,\rm{Id}
$
\\
& $+\;
(t-N/3)
\Bigl(
\frac{p_2 p_3}{
(p_2+p_4)(p_3+p_4)}
-
\frac{p_3p_4(p_1+p_2+p_3+p_4)}{(p_1+p_3)(p_3+p_4)(p_2+p_4)}
\Bigr)$
%(t-N/3)\,\rm{Id}.
\end{tabular}

\end{enumerate}
In the above tables $r=N-s-t$ and $\rho=N-\sigma -\tau$.
\end{theorem}
\noindent {\it Proof:} 
(i) By
(\ref{lem:bilanti-intro}) and since $\varphi^\dagger =\varphi$,
\begin{eqnarray}
\langle \varphi . x^ry^sz^t,
\tilde x^\rho
\tilde y^\sigma
\tilde z^\tau
\rangle
= 
\langle  x^ry^sz^t,
\varphi .
\tilde x^\rho
\tilde y^\sigma
\tilde z^\tau
\rangle.
\label{eq:seven1}
\end{eqnarray}
The left-hand side of
(\ref{eq:seven1}) is equal to
$N! \nu^N$ times
$(s-N/3)P(s,t,\sigma,\tau)$.
To see this,
first observe
$\varphi . x^r
 y^s
 z^t
=
(s-N/3)
 x^r
 y^s
 z^t$
from the tables below
(\ref{eq:Anbasis1-intro}), and then use
Theorem
\ref{thm:pcosines-intro}.
The right-hand side of
(\ref{eq:seven1}) is equal to
$N! \nu^N$ times the weighted sum
in the theorem statement.
To obtain this fact, reduce the
right-hand side of
(\ref{eq:seven1}) using the following three steps:
(a) eliminate $\varphi$ using the long formula in Section 2;
(b) evaluate the result using the tables below
(\ref{eq:Anbasis2-intro});
(c) apply Theorem
\ref{thm:pcosines-intro}.
The result follows from the above comments.
\\
\noindent (ii)--(iv) Similar to the proof of (i) above.
\hfill $\Box$ \\

\noindent
We interpret Theorem
\ref{thm:7termrec} as follows. Parts (i), (ii) 
indicate that
the Rahman polynomials are common eigenvectors
for a pair of difference operators,
while parts (iii), (iv) indicate that 
the dual Rahman polynomials are common eigenvectors
for analogous difference operators. 
 Following \cite{DG} 
 we can say that the 
 Rahman polynomials solve a bispectral problem.
This feature of the Rahman polynomials was hinted at earlier by 
Gr\"unbaum
\cite{grun}.

\section{Directions for future research}\label{se6}

Motivated by the results in this paper
we now pose some problems.
To state the problems we adopt 
a general point of
view.
Let $\F$ denote a field and let
$V$ denote a vector space over $\F$ with finite positive 
dimension.
Let ${\rm End}(V)$ denote the $\K$-algebra consisting
of the $\K$-linear maps from $V$ to $V$.
Fix integers 
$M\geq 1$
and
$N\geq 0$.
Let $\mathbb I=\mathbb I (M,N)$ denote the
set consisting of the $(M+1)$-tuples of nonnegative integers whose sum is
$N$.
A pair of elements
$(r_0,r_1, \ldots, r_M)$ and
$(r'_0,r'_1, \ldots, r'_M)$ in
$\mathbb I$ will be called {\it adjacent} whenever
$(r_0-r'_0,r_1-r'_1, \ldots, r_M-r'_M)$ 
is a permutation of
$(1,-1,0,0,\ldots, 0)$.

\begin{problem}
\label{prob}
\rm
Find all the pairs
 $H$, $\tilde H$ 
that satisfy the following  conditions.
\begin{enumerate}
\item $H$ is an $M$-dimensional subspace of 
${\rm End}(V)$ 
whose elements are diagonalizable
and mutually commute.
\item $\tilde H$ is an $M$-dimensional subspace of 
${\rm End}(V)$ 
whose elements are diagonalizable
and mutually commute.
\item There exists a bijection
$\lambda \mapsto V_\lambda$
from $\mathbb I$ to the set of common eigenspaces of
$H$ such that 
for all $\lambda \in \mathbb I$,
\begin{eqnarray*}
\tilde H V_{\lambda}\subseteq V_\lambda + 
\sum_{\genfrac{}{}{0pt}{} {\mu \in \mathbb I}{\mu \;{\rm adj}\;\lambda}}
V_{\mu}.
\end{eqnarray*}
\item There exists a bijection $\lambda \mapsto \tilde V_\lambda$
from $\mathbb I$ to the set of common eigenspaces of $\tilde H$ such that
for all $\lambda \in \mathbb I$,
\begin{eqnarray*}
 H \tilde V_\lambda \subseteq \tilde V_\lambda +
\sum_{\genfrac{}{}{0pt}{} {\mu \in \mathbb I}{\mu \;{\rm adj}\;\lambda}}
\tilde V_{\mu}.
\end{eqnarray*}
\item
There does not exist a subspace $W$ of $V$ such that
$H W \subseteq W$,
$\tilde H W \subseteq W$,
$W\not=0$,
$W\not=V$.
\item[\rm (vi)] Each of $V_\lambda$,
$\tilde V_\lambda$ has dimension 1 for $\lambda \in \mathbb I$.
\item[\rm (vii)] There exists a nondegenerate symmetric bilinear
form $\langle \,,\,\rangle$ on $V$ such that
both
\begin{eqnarray*}
\langle V_\lambda, V_\mu\rangle &=& 0
\quad {\rm if} \quad \lambda \not = \mu,
\qquad \qquad \lambda, \mu \in \mathbb I,
\label{eq:orth1prob}
\\
\langle \tilde V_\lambda, \tilde V_\mu\rangle &=& 0
\quad {\rm if} \quad \lambda \not = \mu,
\qquad \qquad \lambda, \mu \in \mathbb I.
\label{eq:orth2prob}
\end{eqnarray*}
\end{enumerate}
\end{problem}

\begin{note}
\rm By our results earlier in this paper,
the Rahman polynomials give a solution to
Problem \ref{prob} with $M=2$.
\end{note}

\begin{note}
\rm
For $M=1$ 
a solution to 
Problem 
\ref{prob} is essentially the same thing
as a 
{\it Leonard pair} 
\cite{LS99}.
The Leonard pairs have been studied
extensively and are well understood; see
\cite{madrid} and the references
therein.
The Leonard pairs correspond to a class of orthogonal 
polynomials in one variable.
This class
coincides with the 
terminating branch of the Askey scheme
\cite{KoeSwa}
and consists of the
$q$-Racah, $q$-Hahn, dual $q$-Hahn,
$q$-Krawtchouk,
dual $q$-Krawtchouk, 
quantum 
$q$-Krawtchouk,
affine
$q$-Krawtchouk,
Racah, Hahn, dual Hahn, Krawtchouk,  Bannai/Ito,
and
orphan polynomials
\cite[Section~35]{madrid}.
We expect that for general $M$ 
the solutions to Problem
\ref{prob} correspond to a class of $M$-variable orthogonal
polynomials that resembles the above class.
In particular the $M$-variable polynomials from
\cite[Appendix~A.3]{iliev},
\cite{iliev2},
\cite{tanaka}
are likely members of this class.
\end{note}

\begin{note}
\rm
For $M=1$ the conditions in Problem
\ref{prob} are redundant; indeed
condition (vii) follows from (i)--(vi)
\cite[Lemma~15.2]{madrid}.
\end{note}

\noindent For the ambitious reader
we now pose some harder problems.

\begin{problem}
\label{prob2}
\rm Find all the pairs
$H$, $\tilde H$ that satisfy
conditions (i)--(v) in Problem
\ref{prob}.
\end{problem}

\begin{note}
\rm
For $M=1$ a solution to
Problem \ref{prob2}
is essentially the same thing as a {\it tridiagonal pair}
\cite{TD00};
for a discussion of tridiagonal pairs see
\cite{INT} and the references therein.
The tridiagonal pairs over an algebraically closed field
are classified
\cite[Corollary~18.1]{INT}.
\end{note}

\begin{problem}
\rm 
Above Problem \ref{prob}
we defined an adjacency relation on a 
set $\mathbb I$.
Given its features
a Lie theorist will recognize 
that the solutions to Problems
\ref{prob}, \ref{prob2} belong to the root system $A_M$.
Investigate the analogous objects that belong
to other root systems.
\end{problem}

\section{Acknowledgement}
The second author would like to thank
Hajime Tanaka for several illuminating conversations on
the general subject of this paper, and in particular
for pointing out
reference
\cite{tanaka}. We believe that 
\cite{tanaka} is important and
deserves
to be better known among researchers in special functions.

\bigskip

\noindent Plamen Iliev \hfil\break
\noindent School of Mathematics \hfil\break
\noindent Georgia Institute of Technology \hfil\break
\noindent Atlanta, Georgia  30332-0160 USA \hfil\break
\noindent email: {\tt iliev@math.gatech.edu }\hfil\break

\noindent Paul Terwilliger \hfil\break
\noindent Department of Mathematics \hfil\break
\noindent University of Wisconsin \hfil\break
\noindent 480 Lincoln Drive \hfil\break
\noindent Madison, WI 53706-1388 USA \hfil\break
\noindent email: {\tt terwilli@math.wisc.edu }\hfil\break


\begin{thebibliography}{10}

\bibitem{carter}
R.~Carter.
\newblock {\em Lie algebras of finite and affine type}.
\newblock
Cambridge Studies in Advanced Mathematics 96.
\newblock
Cambridge U. Press. Cambridge, 2005.


\bibitem{DG} J.~J.~Duistermaat and F.~A.~Gr\"unbaum. \newblock Differential 
equations in the spectral parameter. 
\newblock{\em Comm. Math. Phys.} {\bf 103} (1986), 177--240.

\bibitem{iliev}
J.~Geronimo and P.~Iliev.
\newblock Bispectrality of multivariable Racah-Wilson polynomials.
\newblock{\em  Constr. Approx.} {\bf 31} (2010), 417--457;
{\tt arXiv:0705.1469}.


\bibitem{grun}
F.~A.~Gr\"unbaum.
\newblock
The Rahman polynomials are bispectral.
\newblock {\em
 SIGMA Symmetry Integrability Geom. Methods Appl.}
{\bf 3} (2007) Paper 065, 11 pp. 
(electronic).

\bibitem{rahman}
M.R.~Hoare and M.~Rahman.
\newblock A probablistic origin for a new class of bivariate polynomials.
\newblock {\em SIGMA Symmetry Integrability Geom. Methods Appl.}  {\bf 4} 
(2008) Paper 089, 18 pp. (electronic).

\bibitem{iliev2}
P.~Iliev.
\newblock Bispectral commutating  difference operators
for multivariable Askey-Wilson polynomials.
\newblock{\em Trans. Amer. Math. Soc.}, In press.
       {\tt  arXiv:0801.4939}.


%%%%%
%\bibitem{IX} P.~Iliev and Y.~Xu. \newblock Discrete orthogonal polynomials and 
%difference equations of several variables. 
%\newblock{\em Adv. Math.} {\bf 212} (2007), 1--36. {\tt math.CA/0508039}.
%%%%%%%%%%%%%%

\bibitem{INT}
T.~Ito, K.~Nomura, P.~Terwilliger.
\newblock A classification of the sharp tridiagonal pairs.
\newblock{\em Linear Algebra. Appl.}. Submitted.
{\tt arXiv:1001.1812}.

\bibitem{TD00}
T.~Ito, K.~Tanabe, and P.~Terwilliger.
\newblock Some algebra related to ${P}$- and ${Q}$-polynomial association
  schemes,  in:
  \newblock {\em Codes and Association Schemes (Piscataway NJ, 1999)}, Amer.
  Math. Soc., Providence RI, 2001, pp.
       167--192;
       {\tt arXiv:math.CO/0406556}.

%%%%%%%%%%%%%
%\bibitem{Milch} P.~R.~Milch. 
%\newblock A multi-dimensional linear growth birth and 
%death process. \newblock{\em Ann. Math. Statist.} {\bf 39} (1968), 727--754. 
%%%%%%%%%%%%%%%%
\bibitem{janzten}
J.~Jantzen.
\newblock {\em
Lectures on quantum groups}.
\newblock
Graduate Studies in Mathematics, 6.
\newblock
Amer. Math. Soc. Providence, RI, 1996.

\bibitem{KoeSwa}
R.~Koekoek and R.~F.~Swarttouw.
\newblock {\em The Askey scheme of hypergeometric orthogonal
polyomials and its
 $q$-analog}, report 98-17, Delft University of Technology, The
 Netherlands, 1998.
 Available at
\newblock{\tt http://fa.its.tudelft.nl/$\sim$koekoek/askey.html}    

\bibitem{tanaka}
H.~Mizukawa and H.~Tanaka.
\newblock $(n+1,m+1)$-hypergeometric functions associated
to character algebras.
\newblock{\em 
Proc. Amer. Math. Soc.}
{\bf 132} (2004) 2613--2618.


\bibitem{LS99}
P.~Terwilliger.
\newblock Two linear transformations each tridiagonal with respect to an
  eigenbasis of the other.
    \newblock {\em Linear Algebra Appl.}  {\bf 330} (2001) 149--203;
    {\tt arXiv:math.RA/0406555}.


\bibitem{madrid}
P.~Terwilliger.
\newblock
An algebraic approach to the Askey scheme of orthogonal polynomials.
Orthogonal polynomials and special functions,
255--330, Lecture Notes in Math., 1883,
Springer, Berlin, 2006;
{\tt arXiv:math.QA/0408390}.


 \end{thebibliography}
\end{document}